\documentclass{LMCS}

\usepackage{amsmath,amssymb,amsthm}
\usepackage{enumerate}
\usepackage{hyperref}
\usepackage{wasysym}

\usepackage[all]{xy}
\xyoption{2cell}
\xyoption{curve}
\UseTwocells
\input{diagxy}
\usepackage{mathpartir}

\theoremstyle{plain}
\newtheorem*{thmstar}{Theorem}

\newcommand{\A}{A_\bullet}
\newcommand{\uA}[1][]{\underline{A}_{#1}}
\newcommand{\B}{B_\bullet}
\newcommand{\ML}{\mathit{ML_I}}
\newcommand{\MLfrag}{\mathit{ML}^\Id}
\newcommand{\C}{\mathbb{C}}
\newcommand{\bigC}{\mathcal{C}}

\newcommand{\D}{\mathbb{D}}

\renewcommand{\d}{\partial}
\newcommand{\E}{\mathcal{E}}
\newcommand{\F}{\mathcal{F}}
\newcommand{\G}{\mathbb{G}}
\newcommand{\N}{\mathbb{N}}
\renewcommand{\P}{P_{\MLfrag}}
\newcommand{\operadP}{P_{\MLfrag}}

\renewcommand{\r}{{r}}
\renewcommand{\S}{\mathcal{S}}
\newcommand{\T}{\mathcal{T}}
\newcommand{\X}{X_\bullet}
\newcommand{\x}{\vec x}
\newcommand{\uX}[1][]{\underline{X}_{#1}}

\newcommand{\y}{\vec y}
\newcommand{\yon}{\mathbf{y}}

\newcommand{\Alg}{\mathbf{Alg}}
\newcommand{\Cat}{\mathbf{Cat}}

\newcommand{\End}{\mathrm{End}}

\newcommand{\longGSets}{[\mathbb{G}^\op,\mathbf{Sets}]}
\newcommand{\GSets}{\widehat{\mathbb{G}}}
\renewcommand{\lim}{\varprojlim}

\newcommand{\op}{\mathrm{op}}
\newcommand{\Operads}{\mathbf{Operads}}

\newcommand{\Sets}{\mathbf{Sets}}
\newcommand{\Th}{\mathbf{Th}}
\newcommand{\strwCat}{\mathbf{str}\mbox{-}\omega\mbox{-}\mathbf{Cat}}
\newcommand{\wkwCat}{\mathbf{wk}\mbox{-}\omega\mbox{-}\mathbf{Cat}}

\newcommand{\comp}{\textsc{comp}}
\newcommand{\Contr}{\mathsf{Contr}}
\newcommand{\elim}{\textsc{elim}}
\newcommand{\Exch}{\mathsf{Exch}}
\newcommand{\form}{\textsc{form}}
\newcommand{\Id}{\mathrm{Id}}

\newcommand{\idelim}[5]{J_{#1.#2}(#3,#4,#5)}
\newcommand{\intro}{\textsc{intro}}
\newcommand{\Subst}{\mathsf{Subst}}
\newcommand{\src}{\mathsf{src}}
\newcommand{\scterm}{\textsc{term}}
\newcommand{\tgt}{\mathsf{tgt}}
\newcommand{\type}{\mathsf{type}}
\newcommand{\sctype}{\textsc{type}}
\newcommand{\Weak}{\mathsf{Wkg}}
\newcommand{\Vble}{\mathsf{Vble}}

\renewcommand{\equiv}{\simeq}

\newcommand{\iso}{\cong}
\newcommand{\types}{\vdash}

\newcommand{\lscott}{[\![}
\newcommand{\rscott}{]\!]}

\newcommand{\miniqed}{\pushright{$\Diamond$} \penalty-700 \par\addvspace{\medskipamount}}
\newcommand{\doubleqed}{\pushright{$\Diamond$, \qEd} \penalty-700
\par\addvspace{\medskipamount}}

\def\doi{6 (3:24) 2010}
\lmcsheading%
{\doi}
{1--19}
{}
{}
{Dec.~\phantom01, 2009}
{Sep.~17, 2010}
{}

\begin{document}


\title[Weak $\omega$-categories from ITT]{Weak $\omega$-categories from intensional type theory}  

\author[P.~LeF.~Lumsdaine]{Peter LeFanu Lumsdaine}

\address{Department of Mathematics, Carnegie Mellon Univerisity, Pittsburgh, U.S.A.}
\email{plumsdai@andrew.cmu.edu}




\begin{abstract}
\noindent We show that for any type in Martin-L\"o{}f Intensional Type Theory, the terms of that type and its higher identity types form a weak $\omega$-category in the sense of Leinster.  Precisely, we construct a contractible globular operad $\P$ of definable ``composition laws'', and give an action of this operad on the terms of any type and its identity types.
\end{abstract}

\maketitle
\tableofcontents\vfill\eject

\section{Introduction}

\subsection{Overview}

\noindent Starting with the Hofmann-Streicher groupoid model \cite{hofmann-streicher}, higher categories have emerged as a natural approach to the semantics of intensional Martin-L\"of type theory.  In the globular approach to higher categories, a higher category has objects (``0-cells''), arrows (``1-cells'') between objects, 2-cells between 1-cells, and so on, with various composition operations and laws depending on the kind of category in question (strict or weak, $n$- or $\omega$-, $\ldots$).  The paradigm for semantics of type theory is then (very roughly!) that types (or contexts) are thought of as objects $\lscott A \rscott$, terms $x:A \types \tau(x):B$ as arrows $\lscott \tau \rscott \colon \lscott A \rscott \to \lscott B \rscott$, terms of identity type $\rho : \Id_B (\tau, \tau')$ as 2-cells $\lscott \rho \rscott \colon \lscott \tau \rscott \to/=>/ \lscott \tau' \rscott$, terms $\chi : \Id (\rho,\rho')$ as 3-cells, and so on.

This idea has recently been explored by various authors in various directions: see for instance \cite{gambino-garner}, \cite{garner:2-d-models}, \cite{awodey-warren}.  One such direction is investigating the structures formed by the syntax of type theory.  In particular, it has been suggested that (terms of) any type, considered together with its higher identity types, should carry the structure of a weak $\omega$-category or -groupoid.  We will show that this is indeed the case, using the definition of weak $\omega$-category given by Tom Leinster in \cite{leinster:book} following the approach of Michael Batanin \cite{batanin:natural-environment}.

(Note that for this construction, based on a specific type $A$, the dimensions of cells are always one lower than described above: $0$-cells will be terms $\tau:A$, $1$-cells will be terms $\rho:\Id_A(\tau,\tau')$, and so on.  This comes from the general rule that if $X$, $A$ are objects of an $n$-category $\bigC$, then $\bigC(X,A)$ forms an $(n-1)$-category whose 0-cells are 1-cells of $\bigC$, and so on.)

While writing this paper, I found that Benno van den Berg had independently discovered a similar proof (proposed in 2006 and completed in unpublished work \cite{benno:talk}); a development of this is forthcoming in joint work of van den Berg and Richard Garner \cite{benno-richard}. \\

\noindent {\bf Acknowledgements.} I would like to thank Steve Awodey, Pierre-Louis Curien, Richard Garner, Chris Kapulkin, Benno van den Berg, Michael Warren, and the anonymous referees for helpful conversations, comments, and support in preparing this paper.

\subsection{Outline of the construction} \label{subsec:outline}

(We assume throughout some general familiarity with the concepts of higher category theory, but not with the particular definition of weak $\omega$-category used, which we will recall in detail in later sections, and similarly for the type theory.)

In the globular approach, an $\omega$-category $\bigC$ has a set $C_n$ of ``$n$-cells'' for each $n > 0$.  The $0$- and $1$-cells correspond to the objects and arrows of an ordinary category: each arrow $f$ has source and target objects $a = s(f)$, $b = t(f)$.  Similarly, the source and target of a 2-cell $\alpha$ are a parallel pair of 1-cells $f,g: a \two b$, and generally the source and target of an $(n+1)$-cell are a parallel pair of $n$-cells.

\begin{figure}
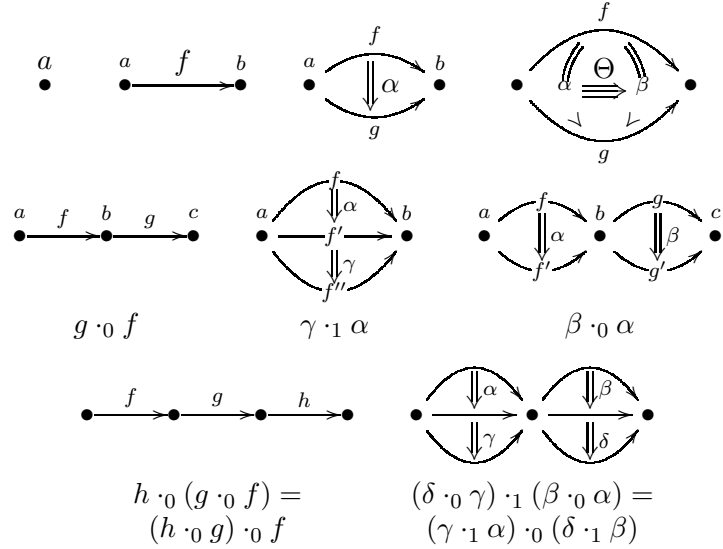

\[
\begin{array}{c}
\begin{array}{cccc}
\ \xy
(0,0)*{\bullet};
(0,80)*{a};
\endxy \quad
&
\ \xy
(0,0)*{\bullet}="a";
(0,80)*{\scriptstyle a};
(400,0)*{\bullet}="b";
(400,80)*{\scriptstyle b};
{\ar "a";"b"};
(200,80)*{f};
\endxy \ 
&
\ \xy
(0,0)*+{\bullet}="a";
(0,80)*{\scriptstyle a};
(450,0)*+{\bullet}="b";
(450,80)*{\scriptstyle b};
{\ar@/^1pc/^{f} "a";"b"};
{\ar@/_1pc/_{g} "a";"b"};
{\ar@{=>} (210,85)*{};(210,-85)*{}};
(280,0)*{\alpha};
\endxy \ 
&
\ \xy
(0,0)*+{\bullet}="a";
(600,0)*+{\bullet}="b";
{\ar@/^1.75pc/^{f} "a";"b"};
{\ar@/_1.75pc/_{g} "a";"b"};
{\ar@2{->}@/_0.5pc/|{\alpha} (220,140);(220,-140)} ;
{\ar@2{->}@/^0.5pc/|{\beta} (380,140);(380,-140)} ;
{\ar@3{->} (225,-20);(375,-20)};
(300,60)*{\Theta};
\endxy \ 
\end{array} \\
\begin{array}{ccc}
\ \xy(0,0)*{\bullet}="a";
(0,80)*{\scriptstyle a};
(300,0)*{\bullet}="b";
(300,80)*{\scriptstyle b};
(600,0)*{\bullet}="c";
(600,80)*{\scriptstyle c};
{\ar^f "a";"b"};
{\ar^g "b";"c"};
\endxy \ 
&
\ \xy
(0,0)*+{\bullet}="a";
(0,80)*{\scriptstyle a};
(500,0)*+{\bullet}="b";
(500,80)*{\scriptstyle b};
{\ar@/^1.75pc/|f "a";"b"};
{\ar|{f'} "a";"b"};
{\ar@/_1.75pc/|{f''} "a";"b"};
{\ar@{=>}^{\alpha} (250,160)*{};(250,50)*{}} ;
{\ar@{=>}^{\gamma} (250,-50)*{};(250,-160)*{}} ;
(0,-250)*{\ };
\endxy \ 
&
\ \xy
(0,0)*+{\bullet}="a";
(0,80)*{\scriptstyle a};
(400,0)*+{\bullet}="b";
(400,80)*{\scriptstyle b};
(800,0)*+{\bullet}="c";
(800,80)*{\scriptstyle c};
{\ar@/^1.1pc/|f "a";"b"};
{\ar@/_1.1pc/|{f'} "a";"b"};
{\ar@/^1.1pc/|g "b";"c"};
{\ar@/_1.1pc/|{g'} "b";"c"};
{\ar@{=>}^{\alpha} (200,80)*{};(200,-80)*{}} ;
{\ar@{=>}^{\beta} (600,80)*{};(600,-80)*{}} ;
\endxy \ \\
g \cdot_0 f &
\gamma \cdot_1 \alpha &
\beta \cdot_0 \alpha
\end{array}
\\
\begin{array}{cc}
\ \xy(0,0)*{\bullet}="a";
(300,0)*{\bullet}="b";
(600,0)*{\bullet}="c";
(900,0)*{\bullet}="d";
{\ar^f "a";"b"};
{\ar^g "b";"c"};
{\ar^h "c";"d"};
\endxy \ &
\ \xy
(0,0)*+{\bullet}="a";
(400,0)*+{\bullet}="b";
{\ar@/^1.5pc/ "a";"b"};
{\ar "a";"b"};
{\ar@/_1.5pc/ "a";"b"};
{\ar@{=>}^{\alpha} (200,150)*{};(200,25)*{}} ;
{\ar@{=>}^{\gamma} (200,-25)*{};(200,-150)*{}} ;
(800,0)*+{\bullet}="c";
{\ar@/^1.5pc/ "b";"c"};
{\ar "b";"c"};
{\ar@/_1.5pc/ "b";"c"};
{\ar@{=>}^{\beta} (600,150)*{};(600,25)*{}} ;
{\ar@{=>}^{\delta} (600,-25)*{};(600,-150)*{}};
(0,250)*{\ };
(0,-220)*{\ };
\endxy \ \\
\begin{array}{c} h \cdot_0 (g \cdot_0 f) =  \\ (h \cdot_0 g) \cdot_0 f \end{array} &
\begin{array}{c}(\delta \cdot_0 \gamma) \cdot_1 (\beta \cdot_0 \alpha) = \\
(\gamma \cdot_1 \alpha) \cdot_0 (\delta \cdot_1 \beta)\end{array}
\end{array}
\end{array}
\]
\caption{Some cells, composites, and associativities in a strict $\omega$-category \label{figure:assoc-laws}} 
\end{figure}

Cells of each dimension can be composed along a common boundary in any lower dimension, and in a \emph{strict} $\omega$-category, the composition satisfies various associativity, unit, and interchange laws, captured by the generalised associativity law: each labelled pasting diagram has a unique composite. (See illustrations in Fig.\ \ref{figure:assoc-laws}).

\noindent In a weak $\omega$-category, we do not expect strict associativity, so may have multiple composites for a given pasting diagram, but we do demand that these composites agree up to cells of the next dimension (``up to homotopy''), and that these associativity cells satisfy certain coherence laws of their own, again up to cells of higher dimension, and so on.

This is exactly the situation we find in intensional type theory.  For instance, even in constructing a term witnessing the transitivity of identity---that is, a composition law for the pasting diagram $(\xymatrix{ \bullet \ar[r] & \bullet \ar[r] & \bullet })$, or explicitly a term $c$ such that 
\[x,y,z:X,\ p:\Id(x,y),\ q:\Id(y,z) \types c(q,p): \Id(x,z)\]
---one finds that there is no single canonical candidate: most obvious are the two equally natural terms $c_l$, $c_r$ obtained by applying ($\Id$-\elim) to $p$ or to $q$ respectively.  These are not definitionally equal, but are propositionally equal, i.e.\ equal up to a 2-cell: there is a term $e$ with
\[x,y,z:X,\ p:\Id(x,y),\ q:\Id(y,z) \types e(q,p): \Id(c_l(q,p),c_r(q,p)).\]

In Leinster's definition \cite{leinster:book}, a system of composition laws of this sort is wrapped up in the algebraic structure of a \emph{globular operad with contraction}, and a weak $\omega$-category is given by a globular set equipped with an \emph{action} of such an operad.  We generalise this slightly, to define an \emph{internal weak $\omega$-category} in any suitable category $\C$.

Accordingly, we would like to find an operad-with-contraction $\operadP$ of all such type-theoretically definable composition laws, acting on terms of any type and its identity types.  In fact, rather than using the full type theory for this, it is more convenient to consider the composition laws definable using just the $\Id$- rules, hence also obtaining the construction for a wider class of theories.

The heart of the paper is Sect.\ \ref{sec:the-operad}, where we formalise this idea.  We consider $\MLfrag[X]$, the fragment of intensional Martin-L\"of type theory generated just by the structural and $\Id$-rules plus a single generic base type $X$. The operad $\operadP$ of definable composition therein laws may then be formally constructed as an endomorphism operad in its syntactic category $\C(\MLfrag [X])$; and by some analysis of the fragment $\MLfrag[X]$, we show that $\operadP$ is contractible.

Since $X$ is generic, $\operadP$ acts on all other types, giving our main theorem:

\begin{thmstar}Let $\T$ be any type theory extending $\MLfrag$, and $A$ any type of $\T$.  Then the system of types $(A, Id_A, Id_{Id_A}, \ldots)$ is equipped naturally with a $\operadP$-action, and hence with the structure of an internal weak $\omega$-category in $\C(\T)$.
\end{thmstar}

To prepare for this, we first lay out in Sect.\ \ref{sec:type-theory-background} our presentation of the type theory $\MLfrag$, and in Sect.\ \ref{sec:operads-background} the relevant background on globular operads and their algebras.

\section{Type-theoretic setting}\label{sec:type-theory-background}

\subsection{The type theories \texorpdfstring{$\MLfrag$}{ML\_Id}, \texorpdfstring{$\MLfrag[X]$}{ML\_Id[X]}}

\noindent Our main theories of interest are the various versions of Intensional Martin-L\"of Type Theory, usually given with identity types ($\Id$-types), dependent sums and products ($\sum$- and $\prod$-types), units ($\mathbf{1}$-types), and possibly more base types (natural numbers, Booleans$\ldots$).  To cover all these in the main theorem, and for a self-contained presentation, we will work throughout this paper in the fragment $\MLfrag$ with only $\Id$-types, and construct our operad from this.

Some care is thus required in our choice of presentation; presentations which are equivalent in the presence of $\sum$- or $\prod$-types may not be so in their absence.  The presentation we use is taken, up to notation, from that of Jacobs \cite{jacobs:categorical-logic}; we list in Table \ref{table:all-rules} the rules assumed, referring to \cite{jacobs:categorical-logic} for their statements, except for the $\Id$-rules, given in full in Table \ref{table:Id-rules}.  A few more of the rules will later be given explicitly, as their precise statements are required.

\begin{table}[hbp]
\begin{center}\begin{tabular}{|@{\ }c@{\ \ \ }c@{\ }|}
\hline
\multicolumn{2}{|c|}{Basic judgement forms} \\
\hline
$ \Gamma \types A\ \type $ & $ \Gamma \types A = B\ \type $ \\
$ \Gamma \types a:A $ & $ \Gamma \types a = b : A $ \\
\hline
\end{tabular} \\
\begin{tabular}{c@{\ \ \ \ }c}
\\ 
\begin{tabular}{|c|}
\hline
Structural groups \\
\hline
Variables ($\Vble$) \\
Substitution ($\Subst$) \\
Weakening ($\Weak$) \\
Exchange ($\Exch$) \\
Equality ($=$) \\
\hline
\end{tabular} &
\begin{tabular}{|c|}
\hline
$\Id$-rules \\
\hline
$\Id$-$\form$ \\
$\Id$-$\intro$ \\
$\Id$-$\elim$ \\
$\Id$-$\comp$ (``$\beta$'' in \cite{jacobs:categorical-logic})\\
compatibility with \\
substitution and $=$ \\
\hline
\end{tabular}
\end{tabular}\end{center}
\caption{The type theory $\MLfrag$} \label{table:all-rules}
\end{table}

\noindent The only features perhaps needing comment are the explicit inclusion of exchange rules, and of the extra dependent context $\Delta$ in the $\Id$-rules; these are each natural rules, but often omitted since they are derivable in the presence of $\Pi$-types (as discussed on e.g.\ p.587 of \cite{jacobs:categorical-logic}).

Note that from $\Exch$ and this $\Id$-\elim\ rule, we can derive a still slightly more general elimination rule $\Id$-$\elim^+$, as $\Id$-\elim\ but with context
\[\Gamma,\ x:A,\ \Delta,\ y:A,\ \Delta',\ p:\Id_A(x,y),\ \Delta''.\]

To simplify notation when referring to iterated identity types, we introduce the notation (following Warren \cite{warren:thesis}) $\uA[n]$ for the $n$th iterated identity type of a type $A$; that is, if $\Gamma \types A\ \type$, then $\Gamma \types \uA[0] := A\ \type$, and inductively
\begin{samepage}
\[\Gamma,\ x_0,y_0:\uA[0],\ x_1,y_1:\uA[1](x_0,y_0),\ \ldots,\
x_n,y_n:\uA[n](x_0,y_0;\ldots ;x_{n-1},y_{n-1}) \qquad \quad
\]
\[ \qquad \qquad \qquad \qquad \qquad \types \uA[n+1](x_0,y_0; \ldots
;x_n, y_n) := \Id_{\uA[n](x_0,\ldots)}(x_n,y_n)\ \type .
\]
\end{samepage}%
We will often omit the superscripts on these when unambiguous.  As usual, we will also be inconsistent in suppression of free variables, writing usually e.g.\ $\y : \Gamma \types A(\y)\ \type$ for clarity in simple cases, but sometimes $\Gamma \types A\ \type$ to avoid unmanageable proliferations of variables.

Finally, for a finite partial order $I = \{i_1 < \ldots < i_n\}$, we will write $\bigwedge_{i \in I} x_i:A_i$ (or just $\bigwedge_{i \in I} A_i$) to denote the context $x_{i_1}:A_{i_1}, \ldots, x_{i_n}:A_{i_n}$.

\begin{table} 
$$\inferrule*[right=$\Id$-\form]{\Gamma \types A\ \type}{\Gamma,\ x,y:A \types \Id_A(x,y)\ \type}$$
$$\inferrule*[right=$\Id$-\intro]{\Gamma \types A\ \type}{\Gamma,\ x:A \types r(x):\Id_A(x,x)}$$
$$\inferrule*[right=$\Id$-\elim]{\Gamma,\ x,y:A,\ p:\Id_A(x,y),\ \Delta(x,y,p) \types C(x,y,p)\ \type \\ \Gamma,\ z:A,\ \Delta(z,z,r(z)) \types d(z):C(z,z,r(z))}{\Gamma,\ x,y:A,\ p:\Id_A(x,y),\ \Delta(x,y,p) \types \idelim{z}{d}{x}{y}{p} : C(x,y,p)}$$
$$\inferrule*[right=$\Id$-\comp]{\textit{(premises as for $\Id$-$\elim$)}}{\Gamma,\ x:A,\ \Delta(x,x,r(x)) \types \idelim{z}{d}{x}{x}{r(x)} = d(x) : C(x,x,r(x))}$$
\caption{Rules for $\Id$-types} \label{table:Id-rules}
\end{table}

\subsection{Translations and syntactic categories}

For reference on this section (including proofs not given here), see Cartmell \cite{cartmell:generalised-algebraic-theories} and Jacobs \cite{jacobs:categorical-logic}.

From here on, we will consider type theories extending $\MLfrag$; formally, by a \emph{type theory} we will mean a \emph{generalised algebraic theory} in the sense of Cartmell \cite{cartmell:generalised-algebraic-theories}, together with an interpretation of the $\Id$-rules in $\T$.

Recall that a \emph{translation} $F$ from such a type theory $\T$ into a type theory $\S$ consists of suitable mappings of types, terms, and derivable judgements, taking each judgement $\Gamma \types A\ \type$ in $\T$ to a judgement $F(\Gamma) \types F(A)\ \type$ in $\S$, and so on, preserving $\Id$-types and their term-constructors, considered up to definitional equality.  (In other words, it is a morphism of generalised algebraic theories, preserving the interpretation of the $\Id$-rules.)

Given $\T$, we write $\T[X]$ for the result of adjoining to $\T$ a fresh base type $X$
\[\inferrule*[right=$X$-\form]{\ }{ \types X\ \type}\]
with no term formation rules.  For any $\S$, a translation $F \colon \T[X] \to \S$ then consists of a translation $F'\colon \T \to \S$ together with a closed type of $\S$. Stating this universal property precisely, in the particular case that we will need:

\begin{prop} \label{prop:universal-property}If $\S$ is any type theory extending $\MLfrag$, and $A$ any closed type of $\S$, then there is a unique translation $F_{S,A} \colon \MLfrag[X] \to \S$ preserving $\Id$-types and their term-constructors and with $F_{S,A}(X) = A$.
\end{prop}

For any type theory $\T$, there is a \emph{syntactic category} $\C(\T)$, having as objects the closed contexts $\Gamma$ of $\T$, and as arrows $f\colon \Gamma \to \Delta$ suitable strings of terms in context $\Gamma$ (\emph{context maps}), all up to definitional equality.  Moreover, a translation $F \colon \T \to \S$ induces a functor $\C(F)\colon\C(\T) \to \C(\S)$; in other words, we have a functor $\C(-)\colon\Th \to \Cat$.

Context maps are sometimes known as \emph{substitutions}, since substitution along them is an important derived rule: if $\y : \Gamma \types A(\y)\ \type$ and $f: \x:\Gamma \to \Delta$ is a context map, then $\x : \Delta \types A(f(\x)\ \type$.  When suppressing free variables, we will write this as $\Delta \types f^*(A)\ \type$.

We will need a simple proposition on limits in syntactic categories:

\begin{prop} \label{prop:dependent-projections-give-limits}
Suppose $\Gamma = \bigwedge_{i \in I} x_i:A_i$ is a context in $\T$, and $\F \subseteq \mathcal{P}(I)$ a set of subsets of $I$, closed under binary intersection and with $\bigcup \F = I$, such that for each $J \in \F$, $\Gamma_J = \bigwedge_{i \in J} x_i : A_i$ is also a well-formed context.

 Then the $\Gamma_J\!$'s and dependent projections between them give a diagram 
\[\Gamma_{-} \colon (\F, \subseteq)^\op \to \C(\T),\]
and the dependent projections $\Gamma \to \Gamma_J$ express $\Gamma$ as its limit:
\[\Gamma = \lim \! {}_{J \in \F}\ \Gamma_J .\]
Moreover, for a translation $F\colon \T \to \T'$, the functor $\C(F)$ preserves such limits. \qed
\end{prop} 

Here \emph{dependent projections} are the obvious context morphisms from a context to any well-formed subcontext, constructed from the $\Vble$, $\Weak$ and $\Exch$ rules.

A familiar special case asserts that if $\Gamma \types A\ \type$ and $\Gamma \types B\ \type$, then the following square of projections is a pullback:
\[\xymatrix{ \Gamma, x:A, y:B \ar[r] \ar[d] & \Gamma ,y:B \ar[d] \\
\Gamma , x: A \ar[r] & \Gamma } 
\]
The proof of the general proposition is essentially the same.

To relativise the constructions of this section to dependent types and contexts over a (closed) context $\Gamma = \bigwedge_{0 \leq i < n} x_i:A_i$ of $\T$, we can consider the \emph{slice type theory} $\T/\Gamma$, given by adjoining to $\T$ a ``generic term of type $\Gamma$'', i.e.\ $n$ new constant symbols $c_i$ and axioms $\types c_i : A_i(c_0,\ldots,c_{i-1})$.  Closed types (resp.\ terms, contexts) of $\T/\Gamma$ then correspond precisely to types (terms, contexts) of $\T$ in context $\Gamma$.

\section{Globular operads and weak \texorpdfstring{$\omega$}{omega}-categories}\label{sec:operads-background}

\noindent As described in the introduction, we want to describe ``the globular operad of composition laws''.  Accordingly, we recall briefly in this section what a globular operad is, and how it formalises the intuition of a set of composition laws for pasting diagrams with structure specifying how these laws themselves compose.  For a slightly (resp.\ much) fuller treatment, and background on strict higher categories, see Leinster \cite{leinster:survey} (resp.\ \cite{leinster:book}).

\subsection{Globular sets and operads}

A \emph{globular set} $\A$ is a presheaf on the category $\G$ generated by arrows
\[0 \two^{s_0}_{t_0} 1 \two^{s_1}_{t_1} 2 \two \ldots \]
subject to the equations $ss = ts$, $st = tt$ (omitting subscripts on the arrows, as usual).  We thus have the category $\GSets := \longGSets$ of globular sets and natural transformations between them.  More generally, a \emph{globular object} in a category $\C$ is a functor $\A \colon  \G^\op \to \C$.  

Explicitly, a globular set $\A$ has a set $A_n$ of ``$n$-cells'' for each $n \in \N$, and each $(n+1)$-cell $x$ has parallel source and target $n$-cells $s(x)$, $t(x)$, as illustrated in the first line of Fig.~\ref{figure:assoc-laws}.  (Cells $x,y$ of dimension $>0$ are \emph{parallel} if $s(x) = s(y)$ and $t(x) = t(y)$; all $0$-cells are considered parallel.)  For parallel $x,y \in A_n$, we write $A(x,y) := \{ z \in A_{n+1}\ |\ s(z) = x, t(z) = y\}$, the set of $(n+1)$-cells from $x$ to $y$.

Our notation will vary: we will typically call globular objects $\A, \B, \ldots$ when emphasising the point of view of the category $\C$, or $A, B, \ldots$ when working more in $[\G^\op, \C]$.  

\begin{exa} \label{ex:pi-omega}For any topological space $X$, there is a globular set $\Pi_\omega(X)$ in which 0-cells are points of $X$, 1-cells are paths between points, 2-cells are homotopies between paths keeping endpoints fixed, and in general, $n$-cells are suitable maps $H\colon [0,1]^n \to X$, viewed as homotopies between $(n-1)$-cells.
\end{exa}

\begin{exa}For any type $A$ in a type theory $\T$, the contexts 
\[x_0,y_0:A,\ x_1,y_1:\uA[1](x_0,y_0)\ldots,\ z:\uA[n](x_0,y_0\ldots,
x_{n-1},y_{n-1}),
\]
along with their dependent projections, form a globular object $\A$ in $\C(\T)$.
\end{exa}

\begin{figure}[htbp]
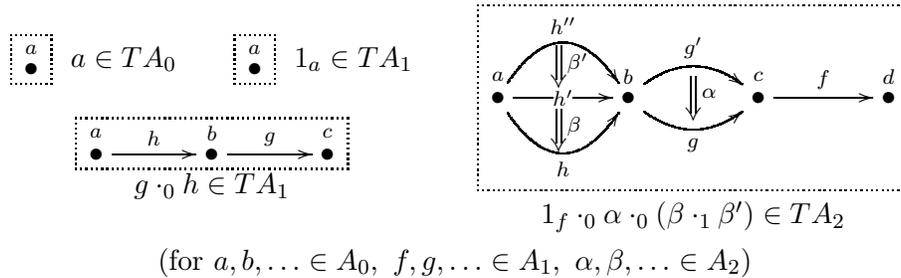

\[\begin{array}{c@{\ \ \ }c}
\begin{array}{c}
\ 
\xy
(0,-35)*{\bullet};
(0,35)*{\scriptstyle a};
(-70,85)*{}="tl";
(70,85)*{}="tr";
(-70,-85)*{}="bl";
(70,-85)*{}="br";
"tl";"tr" **\dir{.};
"tl";"bl" **\dir{.};
"tr";"br" **\dir{.};
"bl";"br" **\dir{.};
\endxy \ \ 
a \in TA_0 \qquad
\xy
(0,-35)*{\bullet};
(0,35)*{\scriptstyle a};
(-70,85)*{}="tl";
(70,85)*{}="tr";
(-70,-85)*{}="bl";
(70,-85)*{}="br";
"tl";"tr" **\dir{.};
"tl";"bl" **\dir{.};
"tr";"br" **\dir{.};
"bl";"br" **\dir{.};
\endxy \ \ 
1_a \in TA_1\ \\ \\
\xy
(0,0)*+{\bullet}="a";
(0,75)*{\scriptstyle a};
(400,0)*+{\bullet}="b";
(400,75)*{\scriptstyle b};
(800,0)*+{\bullet}="c";
(800,75)*{\scriptstyle c};
{\ar^h "a";"b"};
{\ar^g "b";"c"};
(-70,125)="tl";
(870,125)="tr";
(-70,-50)="bl";
(870,-50)="br";
"tl";"tr" **\dir{.};
"tl";"bl" **\dir{.};
"tr";"br" **\dir{.};
"bl";"br" **\dir{.};
\endxy \\
g \cdot_0 h \in TA_1
\end{array}
&
\begin{array}{c}
\xy
(0,0)*+{\bullet}="a";
(0,75)*{\scriptstyle a};
(450,0)*+{\bullet}="b";
(450,75)*{\scriptstyle b};
(900,0)*+{\bullet}="c";
(900,75)*{\scriptstyle c};
(1350,0)*+{\bullet}="d";
(1350,75)*{\scriptstyle d};
{\ar@/^1.75pc/^{h''} "a";"b"};
{\ar|{h'} "a";"b"};
{\ar@/_1.75pc/_h "a";"b"};
{\ar@/^1pc/^{g'} "b";"c"};
{\ar@/_1pc/_g "b";"c"};
{\ar^f "c";"d"};
{\ar@{=>}^<<<{\beta'} (210,180);(210,40)};
{\ar@{=>}^<<<{\beta} (210,-40);(210,-180)};
{\ar@{=>}^<<<{\alpha} (675,75);(675,-75)};
(-70,320)="tl";
(1420,320)="tr";
(-70,-320)="bl";
(1420,-320)="br";
"tl";"tr" **\dir{.};
"tl";"bl" **\dir{.};
"tr";"br" **\dir{.};
"bl";"br" **\dir{.};
(675,-350)*{\ }="strut";
\endxy \\
1_f \cdot_0 \alpha \cdot_0 (\beta \cdot_1 \beta') \in TA_2
\end{array}
\end{array}
\]
\[(\mbox{for }a,b,\ldots \in A_0,\ f,g,\ldots \in A_1,\ \alpha,\beta,\ldots \in A_2)\]
\caption{Some labelled pasting diagrams, elements of a free strict $\omega$-category $T\A$. \label{figure:pastings} } 
\end{figure}

Any strict $\omega$-category (as sketched in the introduction) has an evident underlying globular set, and in fact there is an adjunction (moreover monadic) $F: \GSets \two/->`<-/ \strwCat : U$, giving rise to the ``free strict $\omega$-category'' monad $(T,\mu,\eta)$ on $\GSets$.  Cells of $T\A$ are free (strictly associative) pastings-together of cells from $\A$, including degenerate pastings from the identity cells of $F(\A)$ (as shown in figure \ref{figure:pastings}).

\noindent In particular, $T1$ (where $1$ denotes the terminal globular set, with just one cell of each dimension) consists informally of pastings of this sort, but without labels on the cells.  This is the crucial globular set of \emph{pasting diagrams}.  A peculiarity of $T1$ is that the source and target of any pasting diagram are equal; for this ambivalent operation we write $\d \pi := s(\pi) = t(\pi)$. 

Every pasting diagram $\pi \in T1_n$ has an associated globular set $\hat{\pi}$---intuitively, the set of cells appearing in $\pi$, as shown in our pictures of pasting diagrams throughout.  We then have maps of globular sets $\hat{s}^k, \hat{t}^k\colon  \widehat{\d^k \pi} \to \hat \pi$, embedding $\widehat{\d^k \pi}$ as the $(n-k)$-dimensional source or target of  $\hat \pi$.

Taking categories of elements then gives categories $\int \pi := \int_\G \hat{\pi}$, with objects the cells of $\hat{\pi}$ and arrows into each cell $c$ from its sources and targets $s^k(c)$, $t^k(c)$, and with a functor $\dim \colon  \int \pi \to \G$ giving the dimension of each cell; $\int \pi$ may be seen as the shape of the canonical diagram of basic cells whose colimit in $\GSets$ gives $\hat{\pi}$.

For more discussion of these and other various ways of looking at a pasting diagram, see Street \cite{street:petit-topos}.

A \emph{globular operad} is a globular set $P$ with maps $a\colon P \to T1$ (``arity''), $e\colon  1 \to P$ (``units''), $m \colon  TP \times_{T1} P \to P$ (``composition''), such that
\[\bfig
\node 1(-250,-150)[1]
\node P(250,-150)[P]
\node T1(0,-650)[T1]
\arrow[1`P;e]
\arrow|l|[1`T1;\eta]
\arrow|r|[P`T1;a]
\node TPxP(1000,0)[TP \times_{T1} P]
\node P'(1750,0)[P]
\node TP(750,-250)[TP]
\node P''(1250,-250)[P]
\node T1'(1000,-500)[T1]
\node TT1(850,-650)[T^2 1]
\node T1''(1300,-800)[T1]
\arrow[TPxP`P';m]
\arrow[TPxP`TP;]
\arrow[TPxP`P'';]
\arrow|a|[TP`T1';T!]
\arrow|a|[P''`T1';a]
\arrow|l|/{@{>}@/_5pt/}/[TP`TT1;Ta]
\arrow|l|/{@{>}@/_5pt/}/[TT1`T1'';\mu]
\arrow|r|[P'`T1'';a]
\efig
\]
commute (i.e.\ $e$ and $m$ are maps over $T1$), satisfying the axioms 
\[m \cdot ( \eta \cdot e \times 1_P) = 1_P = m \cdot (\eta \times e)
\colon  P \to P,
\]
\[m \cdot (\mu \times m) = m \cdot (Tm \times 1_P) : T^2 P \times_{T^2
1} TP \times_{T1} P \to P.
\]
Considering the fibers of $a$, we may view $P$ as a family of sets $P(\pi)$ of ``$\pi$-ary operations'' for each $\pi \in T1$: an element $p$ of $P(\pi)$ is seen as a formal operation symbol, taking $\pi$-shaped labelled pasting diagrams as input and returning $n$-cells as output.  The map $e$ then gives us an $n$-cell ``identity'' operation for each $n$, while $m$ allows us to compose operations appropriately.

For readers not familiar with this definition, it may be helpful to first contemplate the simpler case of \emph{plain operads}, defined by diagrams as above but with $\E = \Sets$ and $T$ the ``free monoid'' monad.   These thus have arities valued in $T1 \iso \N$, and present certain finitary, single-sorted equational theories \cite[2.2]{leinster:book}.  However, ``operad'' from here on will always mean ``globular operad''; we will not deal further with any other kind.

A map $f\colon P \to Q$ of globular operads is a map of underlying globular sets commuting with $a, e$ and $m$.

An \emph{action} of a globular operad $P$ on a globular set $A$ is a composition map $c\colon  TA \times_{T1} P \to A$, satisfying 
\[ c \cdot (\eta \times e) = 1_A \colon  A \to A,\]
\[c \cdot (\mu \times m) = c \cdot (Tc \times 1_P) : T^2 A \times_{T^2
1} TP \times_{T1} P \to A.
\]

Informally, this implements the ``formal operations'' in $P$ as actual composition operations on $A$.  An element of $TA \times_{T1} P$ over some $\pi \in T1_n$ is a $\pi$-shaped diagram $\vec x$ with labels from $A$, together with a $\pi$-ary operation $p$ of $P$; $c$ tells us how to apply $p$ to $\vec x$, yielding a single $n$-cell $c(\vec x, p)$ of $A$.

A \emph{$P$-algebra} is a globular set $A$ together with an action of $P$ on $A$. A map $f\colon A \to B$ of $P$-algebras is a map of globular sets commuting with the $P$-actions.  We denote the resulting category by $P$-$\Alg$.

\begin{exa}The globular set $T1$ is itself trivially an operad (indeed, the terminal one), with $a = 1_{T1}$, i.e.\ $T1(\pi) = 1$ for every $\pi$; a $T1$-algebra is then exactly a strict $\omega$-category.  This fits with our description above of a strict $\omega$-category having a unique composition for each pasting diagram.  
\end{exa}

Weak $\omega$-categories will also be described as algebras for a certain globular operad; to find a suitable operad, we need to specify a little extra structure.

A \emph{contraction} on a map $d\colon  A \to B$ of globular sets is a choice of liftings for fillers of parallel pairs: that is, for each parallel pair $x,x' \in A$ (with the convention that all $0$-cells are parallel), a map $\chi_{x,x'}\colon  B(dx,dx') \to A(x,x')$, such that $d \cdot \chi = 1_B$. A \emph{globular operad with contraction} is a globular operad $P$ with a contraction on the map $a\colon  P \to T1$; this ensures both that enough composition operations exist in $P$, and that the operations will be associative up to cells of the next dimension, themselves satisfying appropriate coherence laws up to yet higher cells, and so on.

 It is shown in \cite{leinster:book} that the category of globular operads with contraction has an initial object $L$; this gives the key definition:

\begin{defi} A \emph{weak $\omega$-category} is an $L$-algebra, where $L$ is the initial operad-with-contraction.
\end{defi}

A map $O \to P$ of operads induces a ``restriction of scalars'' functor $P$-$\Alg \to O$-$\Alg$; so if we have an algebra $A$ for any operad $P$ with contraction, restriction along the unique operad-with-contraction map $L \to P$ endows $A$ with the structure of a weak $\omega$-category.

\begin{exa}
The terminal operad $T1$ has a trivial contraction, giving a canonical functor $\strwCat \to \wkwCat$.
\end{exa}

\begin{exa}
For any space $X$, the set $\Pi_\omega(X)$ of Example \ref{ex:pi-omega} may be naturally made into a weak $\omega$-category, the \emph{fundamental weak $\omega$-groupoid} of $X$. \cite[9.2.7]{leinster:book}
\end{exa}

\subsection{Endomorphism operads and more general actions}

\begin{defi} \label{def:a-pi}
For a globular object $\A$ in a category $\C$, and a pasting diagram $\pi \in T1_n$, we define
\[A_\pi := \textstyle \lim_{c \in \int \pi} A_{\dim c},\]
``the object of diagrams of shape $\pi$ in $\A$'', whenever this limit exists in $\C$.  The maps $\hat{s}^k, \hat{t}^k\colon  \widehat{\d^k pi} \to \hat{\pi}$ induce evident projections $s^k, t^k\colon  A_\pi \to A_{\d^k \pi}$.
\end{defi}

An illustration may be useful here: the definition of $A_\pi$ says, for instance, that if $\pi = (\xymatrix{ \bullet \rtwocell & \bullet \rtwocell & \bullet})$, then
\begin{eqnarray*} A_\pi & := & \lim \left( 
\bfig
\node A0l(0,0)[A_0]
\node A0m(700,0)[A_0]
\node A0r(1400,0)[A_0]
\node A1tl(350,300)[A_1]
\node A1tr(1050,300)[A_1]
\node A1bl(350,-300)[A_1]
\node A1br(1050,-300)[A_1]
\node A2l(350,0)[A_2]
\node A2r(1050,0)[A_2]
\arrow[A1tl`A0l;s]
\arrow[A1tl`A0m;t]
\arrow[A1tr`A0m;s]
\arrow[A1tr`A0r;t]
\arrow[A2l`A1tl;s]
\arrow[A2r`A1tr;s]
\arrow[A2l`A1bl;t]
\arrow[A2r`A1br;t]
\arrow[A1bl`A0l;s]
\arrow[A1bl`A0m;t]
\arrow[A1br`A0m;s]
\arrow[A1br`A0r;t]
\efig
\right) \\
& \iso & A_2 \times_{A_0} A_2,
\end{eqnarray*}
giving the object of 0-composable pairs of 2-cells in $A$.  Similarly, if $\pi$ is the basic $n$-cell, then $A_\pi = A_n$.
 
In the case $\C = \Sets$, the sets $A_\pi$ are precisely the fibers of the map $T! \colon  TA \to T1$, by the description of $T$ as a familially representable functor (\cite[8.1]{leinster:book}).


\begin{prop}\label{prop:endo-operad}
If $\A$ is a globular object in a category $\C$, and the objects $A_\pi$ exist, then there is an operad $\End_\C(\A)$, the \emph{endomorphism globular operad of $\A$}, in which (for $\pi \in T1_n$) an element of $\End_\C(\A)(\pi)$ is a sequence of maps $(\sigma_0, \tau_0;\sigma_1,\tau_1;\ldots ,\tau_{n-1}; \rho)$,
\[\rho\colon A_\pi \to A_n, \qquad \sigma_i, \tau_i\colon
A_{\d^{n-i}\pi} \to A_i,
\]
commuting appropriately with the source and target maps, in the sense that
\[s \cdot \sigma_i = s \cdot \tau_i = \sigma_{i-1} \cdot s, \qquad s
\cdot \rho = \sigma_{n-1} \cdot s,
\]
\[t \cdot \sigma_i = t \cdot \tau_i = \tau_{i-1} \cdot t, \qquad t
\cdot \rho = \tau_{n-1} \cdot t.
\]
Moreover, if $F \colon  \C \to \D$ is a functor preserving the limits $A_\pi$, we can also construct $\End(FA)$, and there is a natural map of operads $\End(A) \to \End(FA)$.
\end{prop}

In other words, an element of $\End(\A)(\pi)$ is a map $\rho$ composing diagrams of shape $\pi$ in $\A$ to basic $n$-cells of $A$, extending maps $(\sigma_i,\tau_i)$ composing their sources and targets in each lower dimension.

Such a diagram of maps may be more abstractly seen as a natural transformation $\vec{\rho}\colon  A_{\leq \pi} \to A_{\leq n}$ between two evident functors $A_{\leq \pi}, A_{\leq n}\colon  (\G/n)^\op \to \C$.

\proof
This construction of the endomorphism operad is a straightforward generalisation of the topological case given in \cite[9.2.7]{leinster:book}.  The proof requires more technical background on globular operads from \cite{leinster:book} than can be recalled here; readers unfamiliar with this are encouraged to ``black-box'' this proof and skip to the last few paragraphs of the section.

Recall from \cite[6.4]{leinster:book} that if $S$ is a cartesian monad on a locally cartesian closed category $\E$, then any object $A$ of $\E$ has an \emph{endomorphism $S$-operad} $\End_S(A)$ given by the exponential in $\E/S1$ of the objects $SA \to S1$, $S1 \times A \to S1$; in the internal language of $\E$ this may be written as the dependent sum of exponentials:
\[\End_S(A) = \textstyle \sum_{\pi : S1} [SA_\pi, A].\]

Now, in the case of $(\GSets,T)$, this gives for any globular set $\A$ an endomorphism globular operad $\End_T(\A)$.  For $\pi \in T1_n$, an operad element $p \in \End_T(\A)(\pi)$ then corresponds to a commutative triangle
\[\bfig
\node yn(-300,0)[\yon(n)]
\node EndA(300,250)[\End_T(\A)]
\node T1(300,-250)[T 1]
\arrow[yn`EndA;p]
\arrow|l|[yn`T1;\pi]
\arrow|r|[EndA`T1;]
\efig
\]
and hence, by the definition of $\End_T(\A)$ as a dependent sum of exponentials, to a map from the pullback
\[\bfig
\node TyApi(0,250)[\yon(n) \times_\pi T(\A)]
\node yn(-250,0)[\yon(n)]
\node TyA(250,0)[T(\A)]
\node T1(0,-250)[T1]
\arrow[TyApi`yn;]
\arrow[TyApi`TyA;]
\arrow[yn`T1;\pi]
\arrow[TyA`T1;]
\efig
\] 
into $\A$.  But this in turn corresponds to a map $\yon(n) \times_\pi T(\A) \to \yon(n) \times \A$ in $\GSets/\yon(n)$ and hence, via the equivalence $\GSets/\yon(n) \equiv [(\G/n)^\op,\Sets]$, to a map $\vec{\rho}\colon  A_{\leq \pi} \to A_{\leq n}$ as described above.

Now, given any category $\C$, consider the category $\E = [\C^\op, \GSets] \iso [(\C \times \G)^\op,\Sets] \iso [\G^\op, \widehat{\C}]$.   Composition with $T$ induces a cartesian monad $T^*$ on $\E$.  Since $\E$ is a presheaf category, it is locally cartesian closed; so any object $Y = Y(-)_\bullet$ of $\E$ has an endomorphism $T^*$-operad $\End_{T^*}(Y(-))$.

Moreover, there is an adjunction 
\[\Delta: \E \two/->`<-/ \GSets : \Gamma,\]
where $\Delta$ is the ``$\C$-constant functor'' functor given by $\Delta(\A)(C)_n = A_n$, and $\Gamma$ is ``$\G$-global sections'': $\Gamma(F)_n = \E(\Delta(\yon(n)),F)$.

Using the familial representability of $T$ and the fact that $\Gamma$ preserves limits, we have a cartesian lax map of cartesian monads $(\Gamma,\kappa)\colon  (\E,T^*) \to (\GSets,T)$, and hence an induced functor $\Gamma_* \colon  T^*$-$\Operads_\E \to T$-$\Operads_{\GSets}$.

\begin{defi}Now, any globular object $\A\colon  \G^\op \to \C$ gives an object $\yon(\A)$ of $\E$; we define
\[\End_\C(\A) := \Gamma_* \End_{T^*}(\yon(\A)).\]
\end{defi}

As in the case $\E = \GSets$, we wish to show that this agrees with the explicit description of $\End_\C(\A)$ given in Proposition \ref{prop:endo-operad}.  Again, an element $p \in \End_\C(\A)(\pi)$ is by definition a triangle in $\GSets$:
\[\bfig
\node yn(-300,0)[\yon(n)]
\node EndA(300,250)[\Gamma_* \End_{T^*}(\yon(\A))]
\node T1(300,-250)[\Gamma_* T^* 1]
\arrow[yn`EndA;p]
\arrow|l|[yn`T1;\pi]
\arrow|r|[EndA`T1;]
\efig
\]
which corresponds, by the adjunction $\Delta \dashv \Gamma$ and the definition of $\End_{T^*}(C)$, to a map
$\Delta(\yon(n)) \times_{\Delta(\pi)} T^*(\yon(\A)) \to \Delta(\yon(n)) \times \yon(\A)$

Since all limits and colimits are pointwise, this corresponds to a
family of maps 
\[\yon(n) \times_{\pi} T (\C(C,\A)) \to \yon(n) \times T^*(\C(C,\A))\]
natural in $C$, so (as before) to a natural family of maps
\[\vec \rho_C\colon  \C(C,\A)_{\leq \pi} \to \C(C,\A)_{\leq n},\]
i.e.\ to a map
\[\vec \rho_{(-)}\colon  \yon(\A)_{\leq \pi} \to \yon(\A)_{\leq n}.\]

But since $\yon$ is full and faithful and preserves all existing limits, if the objects $A_\pi$ exist in $\C$ then this corresponds in turn to a map $\vec \rho\colon  A_{\leq \pi} \to A_{\leq n}$, as desired. \qed


We can now extend the definitions of the previous subsection.  An \emph{action} of an operad $P$ on $\A$ is a map of operads $P \to \End(\A)$. (If $\C = \Sets$ then this agrees with our earlier notions of an action on a globular set, by \cite[6.4]{leinster:book}).  A \emph{$P$-algebra} in a category $\C$ is a globular object in $\C$ together with a $P$-action; an \emph{internal weak $\omega$-category} in $\C$ is an $L$-algebra in $\C$.

Moreover, an action of $P$ on $\A$ induces an action of $P$ on the globular set $\C(Y,\A)$ for any $Y \in \C$, since $\C(Y, - ) \colon  \C \to \Sets$ preserves all limits, and hence we have maps $P \to \End_\C(\A) \to \End_\Sets(\C(Y,\A))$.


\section{The contractible globular operad \texorpdfstring{$\P$}{P\_ML\_Id}} \label{sec:the-operad}

\noindent In this section, we construct the promised operad $\P$ of all definable composition laws; we then show that it is contractible, and describe (in the main theorem) how it acts to give the desired weak-$\omega$-category structures on types.

\subsection{Construction of \texorpdfstring{$\P$}{P\_ML\_Id}}

We saw above that for a type $A$ in a type theory $\T$ extending $\MLfrag$, the contexts
\[x_0,y_0:A,\ x_1,y_1:\uA[1](x_0,y_0)\ldots,\ z:\uA[n](x_0,y_0\ldots,
x_{n-1},y_{n-1}),
\]
and the dependent projections between them form a globular context $\A \colon  \G^\op \to \C(\T)$.  In particular, the generic type $X$ gives a globular context $\X$ in $\C(\MLfrag[X])$.  Using the machinery of the previous section, it is now easy to describe $\P$: it will be $\End_{\C(\MLfrag[X])}(\X)$.  However, since $\C(\T)$ does not have all finite limits in general, to use the description of $\End_{\C(\T)}(\A)$ provided by Proposition \ref{prop:endo-operad} we must construct contexts $\Gamma_\pi$ exhibiting the objects $A_\pi$.

Accordingly, suppose we are given $\pi \in T1_n$, with associated globular set $\hat{\pi}$.  There are various ways of putting a total order on the $i$-cells of $\hat{\pi}$ for each $i \leq n$; pick any such.

(There is in fact a canonical choice of such orderings, using the representation of pasting diagrams as Batanin trees (\cite{batanin:natural-environment}, \cite[8.1]{leinster:book}).  This choice has some good compatibility between the orderings on different pasting diagrams, which will later spare us some use of $\Exch$ rules, so for simplicity we will assume it is the ordering chosen; however, since this is purely cosmetic, we will not go into the details here.)

Then take $\Gamma_\pi$ to be the context
\[\bigwedge_{c \in \hat{\pi}_0} x_c\! :\! A,\ \bigwedge_{c \in
\hat{\pi}_1} x_c\! :\! \uA[1](x_{s(c)},x_{t(c)}),\ \ldots\
\bigwedge_{c \in \hat{\pi}_n} x_c\! :\! \uA[n](x_{s^n(c)},x_{t^n(c)};
\ldots;x_{s(c)},x_{t(c)}).
\]
For instance, $\Gamma_{(\bullet \to \bullet \to \bullet)}$ is the context
\[x,y,z:A,\ p:\Id_A(x,y),\ q:\Id_A(y,z)\]
which we met back in the introduction.

Note that we also have projections $\src,\tgt \colon \Gamma_\pi \to \Gamma_{\d \pi}$.

\begin{lem}The context $\x: \Gamma_\pi$, together with the obvious dependent projections, is the object $(\A)_\pi$ of Definition \ref{def:a-pi}; that is, $\Gamma_\pi = \lim_{c \in \int \pi} A_{\dim c}$.  Moreover, if $F \colon \T \to \S$ is a translation of type theories, then $\C(F) \colon \C(\T) \to \C(\S)$ preserves this limit.
\end{lem}
\begin{proof}
Immediate by Proposition \ref{prop:dependent-projections-give-limits} 
\end{proof}

Thus, by Proposition~\ref{prop:endo-operad}, we have:
\begin{thm}The globular object $\A$ in $\C(\T)$ has an endomorphism operad $\End_{\C(\T)}(\A)$, as described in Proposition \ref{prop:endo-operad}, and if $F \colon \T \to \S$ is a translation of type theories, there is an induced map of operads $\End_{\C(\T)}(\A) \to \End_{\C(\S)}(F\A)$. \qed
\end{thm}

Let us unfold what this operad $P := \End_{\C(\T)}(\A)$ actually looks like.  For $\pi \in T1_n$, an element of $P(\pi)$ consists of a map $\rho \colon \Gamma_\pi\to A_n$ in $\C(\T)$, and for $0 \leq k < n$, maps $\sigma_k \colon \Gamma_{\d^{n-k}(\pi)} \to A_k$ and $\tau_k \colon \Gamma_{\d^{n-k}(\pi)} \to A_n$, commuting with the dependent projections.

So, concretely, an element of $P(\pi)$ (a \emph{composition law for $\pi$}) is a sequence of terms $\vec \rho = ((\sigma_i, \tau_i)_{0 \leq i < n}; \rho)$, such that
\begin{eqnarray*}
\x : \Gamma_{\d^n(\pi)} & \types & \sigma_0(\x) : A \\
\x : \Gamma_{\d^n(\pi)} & \types & \tau_0(\x) : A \\
& \vdots & \\
\x : \Gamma_{\d^{n-k}(\pi)} & \types & \sigma_k(\x): \Id (\sigma_{k-1}(\src\ \x),\tau_{k-1}(\tgt\ \x)),\\
\x : \Gamma_{\d^{n-k}(\pi)} & \types & \tau_k(\x): \Id (\sigma_{k-1}(\src\ \x)),\tau_{k-1}(\tgt\ \x)),\\
& \vdots & \\
\x : \Gamma_\pi & \types & \rho(\x) : \Id(\sigma_{n-1}(\src\ \x),\tau_{n-1}(\tgt\ \x)).
\end{eqnarray*} 

The source of this is then the composition law $(\sigma_0,\tau_0,\ldots, \sigma_{n-1},\tau_{n-1};\sigma_n) \in P(\d(\pi))$, and its target is $(\sigma_0,\tau_0,\ldots, \sigma_{n-1},\tau_{n-1};\tau_n) \in P(\d(\pi))$.

We make no attempt to give a formal syntactic description of composition within this operad, but in specific cases it is ``exactly what you would expect'', and is essentially just substitution.  

\begin{defi} \label{defn:operad-p}As as special case of the above construction, we take 
\[\P := \End_{\C(\MLfrag[X])}(\X),\]
 the operad of all definable composition laws on the generic type. 
\end{defi}

For general $\T,A$, we cannot expect $\End_{\C(\T)}(\A)$ to be contractible: contractibility implies (at least) that any two elements of $\End_{\C(\T)}(\A)(\bullet)$ are connected by an element of $\End_{\C(\T)}(\A)(\bullet \to \bullet)$, or in other words that any two terms $x:A \types \tau(x), \tau'(x) : A$ are propositionally equal, which clearly may fail.  However, in the specific case of $\P$, we do wish to show contractibility, since this is the operad which naturally acts on any type.

What precisely does contractibility mean, here?  For every pasting diagram $\pi$ and every parallel pair of composition laws $\vec \sigma, \vec \tau  \in \P(\d(\pi))$, we need to find some filler $\vec \rho \in \P(\pi)$, with $s(\vec \rho) = \vec \sigma$, $t(\vec \rho) = \vec \tau$.

Given $\pi$, such a parallel pair amounts to terms $(\sigma_i,\tau_i)_{0 \leq i < n}$ as in the definition of a composition law for $\pi$, and a filler is a term $\rho$ completing the definition; that is, we seek to derive a judgment
\[\x : \Gamma_\pi \types \rho (\x) : \Id ( \sigma_{n-1}(\src\ \x),
\tau_{n-1} (\tgt\ \x) ).
\]

Playing with small examples (the reader is strongly encouraged to try this---to derive, for instance, the composition and associativity terms mentioned in the introduction) suggests that we should be able to do this by applying $\Id$-\elim\ (possibly repeatedly, working bottom-up as usual) to the variables of identity types in $\Gamma_\pi$.  $\Id$-\elim\ says that to obtain $\rho$, it's enough to obtain it in the case where one of the variables is of the form $r(-)$, and its source and target variables are equal; and by repeated application, it's enough to obtain $\rho$ in the case where multiple higher cells have had identities plugged in in this way.

Now, since the terms $\sigma_i,\tau_i$ have themselves been built up from just the $\Id$-rules, as we plug $r(-)$ terms into them and identify the lower variables, they should sooner or later collapse by $\Id$-\comp\ to be of the form $r^i(x)$ themselves.  In particular, after applying $\Id$-\elim\ as far as possible, plugging in reflexivity terms for the higher variables and contracting all variables of type $X$ to a single $x:X$, the $\sigma_i, \tau_i$ should \emph{all} reduce to reflexivity terms, and in particular $\sigma_{n-1} = \tau_{n-1} = r^{n-1}(x)$, so we can take the desired filler to be
\[x:X \types r^n(x) : \Id(r^{n-1}(x),r^{n-1}(x)).\]
Below, we formalise this argument.  The crucial lemma is that the context $x:X$ is an initial object in $\C(\MLfrag[X])$: that is, since any context $\Gamma$ in $\MLfrag[X]$ is built up just from $X$ and its higher identity types, there is always a unique way to substitute $x$ and its reflexivity terms $r^i(x)$ for all variables of $\Gamma$, and when we subsitute these in to any context morphism $\sigma \colon \Gamma \to \Gamma'$, the result must again reduce to terms of this form.

\subsection{\texorpdfstring{$X$}{X} is initial in \texorpdfstring{$\MLfrag[X]$}{ML\_Id[X]}} \label{subsec:initiality}

\begin{lem} \label{lemma:initiality} The context $x:X$ is an initial object in $\C(\MLfrag[X])$; that is, for any closed context $\Gamma$ there is a unique context map $\ r^\Gamma \colon (x:X) \to \Gamma$. 
\end{lem}

\begin{rem}This lemma does not generally hold in extensions of $\MLfrag[X]$; in $\ML[X]$, for instance, it is easily seen to be false, since for instance there is no term $x:X \types \tau : \Pi_{y:X} \Id(x,y)$.
\end{rem}

\proof We work by structural induction (as, essentially, we must, since this is a property of the theory $\MLfrag[X]$ which can fail in extensions).

So, given any derivation $\delta$ of a judgement $J$ in $\MLfrag[X]$, we recursively derive various terms and/or judgments, depending on the form of $J$, assuming that we have already done so for all sub-derivations of $\delta$.  The form of the terms and judgements we derive will depend on the form of J as follows:\medskip
\[\begin{array}{|c|c|c|c}
\cline{1-3} \rule[-1ex]{0ex}{3.1ex}
J & \textrm{term} & \textrm{judgement} & \\ 

\cline{1-3}  \rule[-1ex]{0ex}{3.5ex} 
\y:\Gamma \types A(\y)\ \type & r^{\Gamma \,\types\, A}(x) & x:X \types r^{\Gamma \,\types\, A}(x) : A(\r^\Gamma (x)) & \\ 

\cline{1-3}  \rule[-1ex]{0ex}{3.5ex} 
\y:\Gamma \types A(\y) = A'(\y) \ \type & - & x:X \types r^{\Gamma \,\types\, A}(x) = r^{\Gamma \,\types\, A'}(x) : A(\r^\Gamma (x)) & (*) \\

\cline{1-3}  \rule[-1ex]{0ex}{3.5ex} 
\y:\Gamma \types \tau(\y) : A(\y) & - & x:X \types \tau(r^\Gamma (x)) = r^{\Gamma \,\types\, A} (x) : A(r^\Gamma (x)) & (**) \\ 

\cline{1-3}  \rule[-1ex]{0ex}{3.5ex} 
\y:\Gamma \types \tau(\y) = \tau'(\y) : A(\y) & - & - & \\ 

\cline{1-3} \end{array}\]\medskip
Here, for a context $\Gamma\ =\ y_0:A_0, \ldots, y_n : A_n(\y_{< n})$, we write $\r^\Gamma$ for the context map $(x:X) \to \Gamma$ consisting of the terms $r^{\types\,A_0}(x) : A_0$, $r^{A_0\, \types\, A_1}(x) : A_1(r^{\types\,A_0}(x))$, \ldots 

Moreover, applying (*) and (**) above to this definition shows that the maps $\r^\Gamma$ respect definitional equality in $\Gamma$, and are preserved by context maps in that for any $f \colon \Delta \to \Gamma$, we have $f(\r^\Delta (x)) = \r^\Gamma (x)$.

Finally, once the induction is complete, applying this last fact together with the definition $r^{(x:X)}(x) := x$ will show that for any other context map $f \colon (x:X) \to \Gamma$, we have $f(x) = f(r^{(x:X)}) = \r^\Gamma (x)$, and so $\r^\Gamma$ is the unique such map, as originally desired.

(This last step is an instance of the general categorical fact that given an object $X$ in a category $\C$ and natural maps $!_Y \colon X \to Y$ to every other object, such that $!_X = 1_X$, it follows that $X$ is initial.)

As usual, the induction proceeds by cases on the last rule used in the derivation of $J$.  Most cases are routine; we include here $X$-\form\ and $\Weak$-\sctype\ as examples of these, together with the less straightforward cases of the $\Id$-rules and $\Subst$-\sctype .

Our definitions for the $\Subst$-\sctype\ and $\Weak$-\sctype\ cases ensure, as usual, that the terms constructed do not depend on the derivation of the judgement used.  As warned earlier, we will vary for readability between showing dependent variables and leaving them implicit, and hence also between the notations $A(f(\x))$ and $f^*A$ for substitution.\medskip

%
%

\noindent($X$-\form):\ in the easiest case, our derivation consists of just the axiom $X$-form
\[\inferrule*[right=$X$-\form]{\ }{ \types X\ \type}\]
and so defining $r^{\types\,X}(X) := x$, we have $x:X \types x: X\ \type$ as needed. \miniqed

\noindent($\Weak$-\sctype):\ Given a derivation ending
\[\inferrule*[right=$\Weak$-\sctype]{\Gamma \types A\ \type \\ \Gamma
\types B\ \type}{\Gamma,\ y:A \types B\ \type}
\]
we inductively already have 
$x:X \types r^{\Gamma \,\types\, B} : (r^\Gamma)^*B $,
and by the $\Subst$ rules,
$ x:X \types (r^\Gamma)^*B = (r^{\Gamma, y:A})^*B\ \type$,
so by equality rules we conclude
$ x:X \types r^{\Gamma \,\types\, B} : (r^{\Gamma,y:A})^*B $
and hence, by $\Weak$-\scterm{}, can set
\[r^{\Gamma,y:A \,\types\, B} := r^{\Gamma \,\types\, B}.\eqno{\Diamond}\]
\noindent($\Id$-\form):\ Given a derivation ending
\[\inferrule*[right=$\Id$-\form]{\Gamma \types A\ \type}{\Gamma,\
y,y':A \types \Id_A(y,y')\ \type}
\]
we need to find a term
\[x:X \types r^{\Gamma,y,y':A \,\types\, \Id_A(y,y')} :
\Id_{(r^\Gamma)^*A}(r^{\Gamma \,\types\, A},r^{\Gamma,y:A \,\types\,
A})
\]
But $\Gamma, y:A \types A\ \type$ may be derived using weakening, and so by our construction for $\Weak$-\sctype\ above, $r^{\Gamma,y:A \,\types\, A} = r^{\Gamma \,\types\, A}$, so we have
\[x:X \types r(r^{\Gamma \,\types\, A}) :
\Id_{(r^\Gamma)^*A}(r^{\Gamma \,\types\, A},r^{\Gamma,y:A \,\types\,
A})
\]
and so can set
\[r^{\Gamma, y,y':A \,\types\, \Id_A(y,y')} := r(r^{\Gamma \,\types\,
A}). \eqno{\Diamond}
\]
\noindent($\Id$-\intro):\ Now we are given a derivation with last step
\[\inferrule*[right=$\Id$-\intro]{\Gamma \types A\ \type}{\Gamma,\ y:A
\types r(y):\Id_A(y,y)}
\]
and wish to show
\[x:X \types r(r^{\Gamma \,\types\, A}) = r^{\Gamma, y,y:A \,\types\,
\Id_A(y,y)} : (r^\Gamma)^* A.
\]
But by our construction of $r^{\Gamma, y,y':A \,\types\, \Id_A(y,y')}$ above (our $\Id$-\form\ case), and of $r^{\Gamma, y:A \,\types\, \Id_A(y,y)}$ from it (our $\Contr$-\sctype\ case), this is just the definition of $r^{\Gamma, y,y:A \,\types\, \Id_A(y,y)}$. \miniqed \medskip


\noindent($\Id$-\elim):\ Here, we are given a derivation ending
\[\inferrule*[right=$\Id$-\elim]{\Gamma,\ y,y':A,\ p:\Id_A(y,y'),\
\Delta(y,y',p) \types C(y,y',p)\ \type \\ \Gamma,\ z:A,\
\Delta(z,z,r(z)) \types d(z):C(z,z,r(z))}{\Gamma,\ y,y':A,\
p:\Id_A(y,y'),\ \Delta(y,y',p) \types \idelim{z}{d}{y}{y'}{p} :
C(y,y',p)};
\]
for readability, we assume $\Delta$ is empty.  We want to derive the judgement
\begin{eqnarray*} x:X \types (r^{\Gamma, y,y':A, p:\Id(y,y')})^* (\idelim{z}{d}{y}{y'}{p}) & = & r^{\Gamma, y,y':A, p:\Id(y,y') \,\types\, C(y,y',p)} \\
& & \qquad : (r^{\Gamma, y,y':A, p:\Id(y,y')})^*C.
\end{eqnarray*}
Unwrapping the former term, we have (all in context $(x:X)$):
\[\begin{tabular}{llr}
\multicolumn{3}{l}{$\displaystyle (r^{\Gamma, y,y':A, p:\Id(y,y')})^* (\idelim{z}{d}{y}{y'}{p})$} \\
$\quad$ & $\displaystyle = \idelim{z}{(r^\Gamma)^* d}{r^{\Gamma \,\types\, A}}{r^{\Gamma,y:A \,\types\, A}}{r^{\Gamma,y,y':A \,\types\, \Id(y,y')}}$& \\
& $\displaystyle = \idelim{z}{(r^\Gamma)^* d}{r^{\Gamma \,\types\, A}}{r^{\Gamma \,\types\, A}}{r(r^{\Gamma \,\types\, A})}$ & \\
& $\displaystyle = (r^\Gamma)^* d(r^{\Gamma \,\types\, A})$ & (by $\Id$-\comp)\\
& $\displaystyle = (r^{\Gamma, z:A})^* d$ & \\
& $\displaystyle = r^{\Gamma, z:A \,\types\, C(z,z,r(z))}$ & (by induction) \\
& $\displaystyle = r^{\Gamma, y,y':A, p:\Id(y,y') \,\types\, C(y,y',p)}$ \\
\multicolumn{3}{r}{(by the definition of $\displaystyle r^{\Gamma, z:A \,\types\, C(z,z,r(z))}$} \\
\multicolumn{3}{r}{using our $\Weak$-\sctype\ and $\Id$-\elim\ cases.)}
\end{tabular}
\]
If $\Delta$ in the application of $\Id$-\elim\ is non-empty, we have a few more lines, relying inductively on our $\Subst$-rules cases. \miniqed\medskip

\noindent($\Subst$-\sctype):\ For this case we will need one more piece of notation, generalising the context maps $\r^\Gamma$: for a dependent context $\Delta = \bigwedge_i A_i$ over $\Gamma$, we write $r^{\Gamma \,\types\, \Delta} \colon (x:X) \to (r^\Gamma)^* \Delta$ for the map built up from terms $r^{\Gamma, \Delta_{< i} \,\types\, A_{i}}$ in the obvious way.

So, we are given a derivation ending with the rule
\[\inferrule*[right=$\Subst$-\sctype]{\Gamma,\ y:A,\ \Delta \types B\
\type \\ \Gamma \types f:A}{\Gamma,\ f^*\Delta \types f^*B\ \type}
\]
and we wish to derive a judgement
\[x:X \types r^{\Gamma,f^*\Delta \,\types\, f^*B} :
(f^{\Gamma,f^*\Delta})^*(f^*B).
\]
Unfolding the definition of the desired type, we have 
\[\begin{tabular}{ll}
\multicolumn{2}{l}{$\displaystyle (f^{\Gamma,f^*\Delta})^*(f^*B)$} \\
$\quad$ & $\displaystyle = (r^{\Gamma \,\types\, f^*\Delta})^*(r^\Gamma)^*(f^*B)$ \\
 & $\displaystyle = (r^{\Gamma \,\types\, f^*\Delta})^*((r^\Gamma)^*f)^*(r^\Gamma)^*B$ \\
 & $\displaystyle = (r^{\Gamma \,\types\, f^*\Delta})^*(r^{\Gamma \,\types\, A})^*(r^\Gamma)^*B \qquad$ \hfill (by induction) \\
 & $\displaystyle = (r^{\Gamma \,\types\, f^*\Delta})^*(r^{\Gamma,y:A})^*B$ \\
 & $\displaystyle = (r^{\Gamma,y:A \,\types\, \Delta})^*(r^{\Gamma,y:A})^*B \qquad \qquad \quad$ \hfill (by def'n of $r^{\Gamma \,\types\, f^*\Delta}$, i.e.\ by previous\\
 & \hfill applications of \emph{this} case)  \\
 & $\displaystyle = (r^{\Gamma,y:A,\Delta})^*B$ \\
\end{tabular}
\]
so since by induction $\displaystyle x:X \types r^{\Gamma,y:A,\Delta \,\types\, B}:(r^{\Gamma,y:A,\Delta})^*B$, we take $r^{\Gamma,f^*\Delta \,\types\, f^*B} := r^{\Gamma,y:A,\Delta \,\types\, B}$.

The cases for the other structural rules and $X$-\form\ are straightforward, similar to the $\Weak$-\sctype\ case above. \doubleqed

\subsection{Contractibility of \texorpdfstring{$\P$}{P\_ML\_Id}}

We are now ready to show that $\P$ is contractible, arguing along the lines sketched above.

\begin{thm}\label{theorem:p-is-contractible}The operad $\P$ is contractible.
\end{thm}

\begin{proof}As described above, this amounts to the statement: for every $n \in \N$ and pasting diagram $\pi \in T1_n$, and every sequence $(\sigma_i,\tau_i)_{i<n}$ of terms such that
\[\x : \Gamma_{\d^{n-i}(\pi)} \types \sigma_i(\x) : \uX
(\sigma_0(\src^i\ \x),\ldots,\tau_{i-1}(\tgt\ \x))
\]
\[\x : \Gamma_{\d^{n-i}(\pi)} \types \tau_i(\x) : \uX
(\sigma_0(\src^i\ \x),\ldots,\tau_{i-1}(\tgt\ \x))
\]
($i < n$) are derivable in $\MLfrag[X]$, we can find a ``filler'', i.e.\ a term $\rho$ with
\[\x : \Gamma_\pi \types \rho(\x) : \uX (\sigma_0(\src^n\
\x),\ldots,\tau_{n-1}(\tgt\ \x))
\]
We show this by induction on the number of cells in $\pi$.

Suppose $\pi$ has more than one cell.  Then it must have some cells in dimension $> 1$.  Let $k$ be the highest dimension in which $\pi$ has cells, and $c$ be some $k$-cell of $\hat{\pi}$.  Now take $\pi^{-c} \in T1_n$ to be the pasting diagram whose globular set is obtained (up to isomorphism) from that of $\pi$ by removing $c$ and identifying $s(c)$ and $t(c)$.

Now $\Gamma_{\pi^{-c}}$ is exactly (up to renaming of variables, and possibly re-ordering if we do not assume that we chose compatible orderings of the cells of pasting diagrams) the context obtained from $\Gamma_\pi$ by removing the variables $x^k_c$ and $x^{k-1}_{t(c)}$, and replacing any occurrences of the latter in subsequent types by $x^{k-1}_{s(c)}$, and we have a natural context map $h \colon \Gamma_{\pi^{-c}} \to \Gamma_\pi$ given by plugging in $x^{k-1}_{s(c)}$ for $x^{k-1}_{t(c)}$ and $r(x^{k-1}_{s(c)})$ for $x^k_c$; and these are exactly right for
\[\inferrule*{\x:\Gamma_{\pi^{-c}} \types \rho^{-c}(\x) : \uX
(\sigma_0(\src^n\ h(\x)),\ldots,\tau_{n-1}(\tgt\ h(\x)))}{\x :
\Gamma_\pi \types
\idelim{x^{k-1}_{s(c)}}{\rho^{-c}}{x^{k-1}_{s(c)}}{x^{k-1}_{t(c)}}{x^k_c}
: \uX (\sigma_0(\src^n\ \x),\ldots,\tau_{n-1}(\tgt\ \x))}
\]
to be an instance of $\Id$-$\elim^+$.  So to give the desired filler $\rho$, it is enough to give $\rho^{-c}$ with
\[\x:\Gamma_{\pi^{-c}} \types \rho^{-c}(\x) : \uX (\sigma_0(\src^n\
h(\x)),\ldots,\tau_{n-1}(\tgt\ h(\x))).
\]

But now note that
\[\d^{n-i}(\pi^{-c}) = \left\{ \begin{array}{ll} \d^{n-i}(\pi) &
\textrm{for $n-i < k$} \\ (\d^{n-i}(\pi))^{-c} & \textrm{for $n-i \geq
k$} \end{array} \right. ;
\] 
moreover, we can construct context maps
\[h^s_i, h^t_i \colon \Gamma_{\d^{n-i}(\pi^{-c})} \to
\Gamma_{\d^{n-i}(\pi)}
\]
(analogous to $h$ if $i \geq k$, and just the identity otherwise), and these commute with the maps $\src$ and $\tgt$.  So for each $i < n$, we have
\[\x : \Gamma_{\d^{n-i}(\pi^{-c})} \types \sigma_i(h(\x)) : \uX
(\sigma_0(h(\src^i\ \x)),\ldots,\tau_{i-1}(h(\tgt\ \x))),
\]
\[\x : \Gamma_{\d^{n-i}(\pi^{-c})} \types \tau_i(h(\x)) : \uX
(\sigma_0(h(\src^i\ \x)),\ldots,\tau_{i-1}(h(\tgt\ \x))),
\]
i.e.\ the sequence of terms $(h^*(\sigma_i),h^*(\tau_i))_{i<n}$ are a parallel pair for $\pi^{-c}$.  So by induction (since $\pi^{-c}$ has fewer cells than $\pi$), these terms have a filler; but this filler is exactly the desired term $\rho^{-c}$.

Thus it is enough to show the existence of fillers in the case where $\pi$ has just one cell, i.e.\ where $\pi = ( \bullet )$.  But in this case, $\Gamma_\pi = \Gamma_{\d^i(\pi)} = \Gamma_{\d^i(\pi)} = (x:X)$ for each $i < n$, and so by the initiality of $(x:X)$ we must have $\sigma_i(x) = \tau_i(x) = r^i(x)$ for each $i$; so now $\rho := r^n(x)$ gives the filler, and we are done.
\end{proof}

Unwinding this induction, we can see that it exactly formalises the process described at the start of Subsection \ref{subsec:initiality}, of repeatedly plugging in higher reflexivity terms for all variables, knowing that the given composites will themselves eventually compute down to higher reflexivity terms.

Note that Lemma \ref{lemma:initiality} was applied only at the base case of the induction, and only to show that terms $x:X \types \sigma: \Id(r^n(x),r^n(x))$ must be equal to $r^{n+1}(x)$.  A sufficiently strong normalisation result would also imply this, resting on showing that these are the only appropriate normal forms; this could then extend also to the operad $\End_{\ML[X]}(\X)$ of all composition laws of the \emph{full} type theory, which cannot be shown contractible by the present method.  However, working with the fragment $\MLfrag$ seems more economical, showing that $\Id$-types are the only structure required.

\subsection{Types as weak \texorpdfstring{$\omega$}{omega}-categories} \label{subsec:payoff}

Putting the above results together, we obtain our main goal:

\begin{thm}Let $\T$ be any type theory extending the fragment $\MLfrag$, $\Gamma$ any closed context of $\T$, $A$ a dependent type over $\Gamma$.  Then the globular context $\A$ carries the structure of a $\P$-algebra in $\C(\T/\Gamma)$.
\end{thm}

\proof By Proposition \ref{prop:universal-property}, there is a unique translation $F_{\T/\Gamma,A} \colon  \MLfrag[X] \to \T/\Gamma$ taking $X$ to $A$, and hence taking $\X$ to $\A$.  By Proposition \ref{prop:endo-operad}, this induces an action of $\P$ on $\A$, and so, since by Theorem \ref{theorem:p-is-contractible} $\P$ admits a contraction, an action of $L$ (the initial operad-with-contraction) on $\A$, as desired. \qed

\begin{cor}Let $\T$, $\Gamma$, $A$ be as above, and $\Delta$ a dependent context over $\Gamma$.  Then the globular set of terms of types $A$, $\Id_A$, $\Id_{\Id_A}$, $\ldots$ in context $\Gamma, \Delta$ carries the structure of a $\P$-algebra, and hence of a weak $\omega$-category.
\end{cor}

\proof This is just the globular set of $\C(\T/\Gamma)(\Delta,\A)$ of context maps
\[f\colon  \Gamma, \Delta \to \Gamma, x_0,y_0:A, \ldots ,\
x_{n-1},y_{n-1}:\uA(x_0,\ldots,y_{n-2}),\ z:\uA(x_0,\ldots,y_{n-1})
\]
and so inherits a $\P$-action, and hence an $L$-action, from the actions on $\A$.
\qed

\begin{rem}[Functoriality]  The construction of the $\P$-algebra $\C(\T/\Gamma)(\Delta,\A)$ should be covariantly functorial in $\T$, and contravariantly in $\Gamma$ and $\Delta$.  That is, translations $\T \to \T'$ and context maps $\Gamma' \to \Gamma$, $\Delta' \to \Delta$ should induce \emph{strict} maps of $\P$-algebras, composing appropriately.  A proof of this should be fairly straightforward, by an extension of the methods of the current paper; essentially, the missing ingredient is a treatment of maps of internal operad algebras.

More subtly, it should be functorial in $A$, but only to \emph{weak} maps:  a map of types $A \to A'$ should induce \emph{weak} maps of $\P$- or $L$-algebras---that is, weak $\omega$-functors.  This seems an altogether trickier question, due partly but not only to the lack, until fairly recently (\cite{garner:homomorphisms}), of a suitable definition of weak $\omega$-functor.
\end{rem}

\begin{rem}[Comparison with \cite{benno-richard}]  As mentioned in the Introduction, Richard Garner and Benno van den Berg have independently given (\cite{benno-richard}) a proof of essentially the same result.  The core of their approach is the same as that given here: the $\omega$-category action is induced via contractible operads constructed from endomorphism operads of the globular contexts of identity types.  The main differences between construction of the present paper and that of \cite{benno-richard} are, roughly, as follows:

\begin{enumerate}[(1)]
\item Garner and van den Berg use Batanin's presentation \cite{batanin:natural-environment} of globular operads and higher categories, while I have used the later presentation of Leinster \cite{leinster:book}.  This is essentially a superficial difference; the two presentations  are intertranslatable.
\item  Garner and van den Berg work from the categorical structure on syntactic categories given by the identity types, rather than from the identity types in the syntax directly.
\item  Where I have used the single operad $\P$ of definable composition laws on the generic type, Garner and van den Berg use, for each type, a tailor-made operad of composition laws \emph{on that type}, constructed from the endomorphism operad over it in the syntactic category of the particular theory in question.
\item  As remarked after Definition \ref{defn:operad-p}, this entire endomorphism operad will not in general be contractible; consequently, Garner and van den Berg pass to a sub-operad of ``point-preserving'' operations, which is always contractible.  From this point of view, Subsection \ref{subsec:initiality} (the initiality of $(x:X)$ in $\MLfrag[X]$) may be seen as showing that over the generic type in $\MLfrag[X]$, \emph{all} composition laws are point-preserving.
\item  Finally, Garner and van den Berg show moreover that the weak $\omega$-categories produced are in fact weak $\omega$-\emph{groupoids}, according to the criterion of Cheng \cite{cheng:duals-give-inverses}.
\end{enumerate}
\end{rem}

\bibliographystyle{amsalpha}
\bibliography{higher-cats-type-theory-bib}

\providecommand{\bysame}{\leavevmode\hbox to3em{\hrulefill}\thinspace}
\providecommand{\MR}{\relax\ifhmode\unskip\space\fi MR }
\providecommand{\MRhref}[2]{%
  \href{http://www.ams.org/mathscinet-getitem?mr=#1}{#2}
}
\providecommand{\href}[2]{#2}
\begin{thebibliography}{GvdB08}

\bibitem[AW09]{awodey-warren}
Steve Awodey and Michael~A. Warren, \emph{Homotopy theoretic models of identity
  types}, Math. Proc. Cambridge Philos. Soc. \textbf{146} (2009), no.~1,
  45--55. \MR{MR2461866}

\bibitem[Bat98]{batanin:natural-environment}
M.~A. Batanin, \emph{Monoidal globular categories as a natural environment for
  the theory of weak {$n$}-categories}, Adv. Math. \textbf{136} (1998), no.~1,
  39--103.

\bibitem[Car86]{cartmell:generalised-algebraic-theories}
John Cartmell, \emph{Generalised algebraic theories and contextual categories},
  Ann. Pure Appl. Logic \textbf{32} (1986), no.~3, 209--243.

\bibitem[Che07]{cheng:duals-give-inverses}
Eugenia Cheng, \emph{An {$\omega$}-category with all duals is an
  {$\omega$}-groupoid}, Appl. Categ. Structures \textbf{15} (2007), no.~4,
  439--453.

\bibitem[Gar08]{garner:homomorphisms}
Richard Garner, \emph{Homomorphisms of higher categories}, submitted, 2008.

\bibitem[Gar09]{garner:2-d-models}
\bysame, \emph{Two-dimensional models of type theory}, Math. Structures Comput.
  Sci. \textbf{19} (2009), no.~4, 687--736. \MR{MR2525957}

\bibitem[GG08]{gambino-garner}
Nicola Gambino and Richard Garner, \emph{The identity type weak factorisation
  system}, Theoret. Comput. Sci. \textbf{409} (2008), no.~1, 94--109.
  \MR{MR2469279}

\bibitem[GvdB08]{benno-richard}
Richard Garner and Benno van~den Berg, \emph{Types are weak
  {$\omega$}-groupoids}, submitted, 2008.

\bibitem[HS98]{hofmann-streicher}
Martin Hofmann and Thomas Streicher, \emph{The groupoid interpretation of type
  theory}, Twenty-five years of constructive type theory ({V}enice, 1995),
  Oxford Logic Guides, vol.~36, Oxford Univ. Press, New York, 1998,
  pp.~83--111.

\bibitem[Jac99]{jacobs:categorical-logic}
Bart Jacobs, \emph{Categorical logic and type theory}, Studies in Logic and the
  Foundations of Mathematics, vol. 141, North-Holland Publishing Co.,
  Amsterdam, 1999.

\bibitem[Lei02]{leinster:survey}
Tom Leinster, \emph{A survey of definitions of {$n$}-category}, Theory Appl.
  Categ. \textbf{10} (2002), 1--70 (electronic).

\bibitem[Lei04]{leinster:book}
\bysame, \emph{Higher operads, higher categories}, London Mathematical Society
  Lecture Note Series, vol. 298, Cambridge University Press, Cambridge, 2004.

\bibitem[Str00]{street:petit-topos}
Ross Street, \emph{The petit topos of globular sets}, Journal of Pure and
  Applied Algebra \textbf{154} (2000), 299--315.

\bibitem[vdB]{benno:talk}
Benno van~den Berg, \emph{Types as weak {$\omega$}-categories}, Lecture
  delivered in Uppsala, 2006, and unpublished notes.

\bibitem[War08]{warren:thesis}
Michael~A. Warren, \emph{Homotopy theoretic aspects of constructive type
  theory}, Ph.D. thesis, Carnegie Mellon University, 2008.

\end{thebibliography}
\end{document}